\input amstex.tex 
\documentstyle{amsppt} 
\topmatter
\title
K3 surfaces via almost-primes
\endtitle
\author
Keiji Oguiso
\endauthor
\address 
\endaddress

\abstract Based on the result on derived categories on K3 surfaces due to Mukai and Orlov and the result concerning almost-prime numbers due to Iwaniec, we remark the following fact: For any given positive integer $N$, there are $N$ (mutually non-isomorphic) projective complex K3 surfaces such that their 
Picard lattices are not isomorphic but their transcendental lattices are Hodge isometric, or equivalently, their derived categories are mutually equivalent. 
After reviewing finiteness result, we also give an explicit formula for the cardinality of the isomorphism classes of projective K3 surfaces having 
derived categories equivalent to the one of $X$ with Picard number $1$ in 
terms of the 
degree of $X$.  
\endabstract

\leftheadtext{K. Oguiso}
\rightheadtext{K3 surfaces via almost-primes}
\endtopmatter

\document
\head
{\S 1. Introduction}
\endhead 
Our main results are (1.7) - (1.10) and Proposition A in the appendix. 
\flushpar
\vskip 4pt
(1.1) Let $\bold N := \{1, 2, 3, \cdots \}$ be the set of natural numbers. We call $p \in \bold N$ {\it prime} if 
$n \vert p$ and $n \in \bold N$ then $n = 1$ or $p$. In this note it is convenient to consider the number $1$ is also prime and we shall do so. There is a long standing conjecture since Dirichlet concerning primes: 
\proclaim{Conjecture} Set $\Cal N_{1} := \{ n \in \bold N \vert n^{2} + 1\, \text{is prime}\, \}$. Then $\vert \Cal N_{1} \vert = \infty$. \qed \endproclaim 
Contrary to its appearance, this conjecture is very difficult and is completely open till now. The products of two (not necessarily distinct) primes are called {\it almost-primes}. The best result known now is the following very deep Theorem due to H. Iwaniec:
\proclaim{Theorem (1.2) [Iw]} Let $f(x) = ax^{2} + bx + c$ be an irreducible polynomial in $\bold Z[x]$ such that $c \equiv 1 (2)$ and $a > 0$. Then for each such an $f$ the set of almost-primes of the form $f(n)$ ($n \in \bold N$) are infinite. In particular, almost-primes of the form $n^{2} + 1$ are infinite. \qed \endproclaim 
This Theorem concerning quadratically represented almost-primes is one of the most crucial ingredients of this note. 
In particular, we apply one special case of this Theorem, namely, the case $f(x) = 4x^{2} +1$ for our study of K3 surfaces. (The argument also goes through if we take $f(x) = x^{2} + 1$, but then we sometimes need case by case description according to the congruence $\text{mod}\, 4$. This is only the reason why we choose $4x^{2} + 1$ instead of more attractive $x^{2} + 1$.) 
We set 
$$\Cal N_{2} := \{n \in \bold N \vert 4n^{2} + 1\,\, \text{is almost prime}\}.$$ 
By (1.2), $\vert \Cal N_{2} \vert = \infty$. If $n \in \Cal N_{2}$, then we have the factorization $4n^{2} + 1 = pq$. Note that $4n^{2} + 1 \equiv 1\, \text{mod}\, 4$, $p \equiv q \equiv 1\, \text{mod}\, 2$ and $p \not= q$; Indeed, if $p = q$, then $1 = (p-2n)(p+2n)$, a contradiction to $p + 2n \geq 2$. 
\flushpar
\vskip 4pt
(1.3) Throughout this note, a {\it K3 surface} means a $2$-dimensional complex {\it projective} smooth variety which admits nowhere vanishing global 2-form but admits no nonzero global one form. Let $X$ be a K3 surface. Let $(* , **)$ be the symmetric bilinear form on 
$H^{2}(X, \bold Z)$ given by the cup product. Then $(H^{2}(X, \bold Z), (* , **))$ is an even unimodular lattice of signature $(3, 19)$. We denote by $\text{NS}\,(X) \simeq \text{Pic}\, X$ the N\'eron-Severi lattice of $X$ and by $\rho(X)$ the Picard number, i.e. rank of $\text{NS}(X)$. This lattice $\text{NS}\,(X)$ is primitive in 
$H^{2}(X, \bold Z)$ and is of signature $(1, \rho(S) -1)$. 
We call the orthogonal lattice of $\text{NS}\,(X)$ in $H^{2}(X, \bold Z)$ the transcendental lattice and denote it by $T_{X}$, i.e.
$$T_{X} := \{x \in H^{2}(X, \bold Z) \vert (x, y) = 0\, \text{for all}\, 
y \in \text{NS}\,(X)\}.$$ 
$T_{X}$ is also primitive in $H^{2}(X, \bold Z)$ and is of signature 
$(2, 20 - \rho(X))$. Let $\omega_{X}$ be a non-zero global 2-form on $X$. Then, one has the inclusion
$$\bold C \omega_{X} \oplus \bold C \overline{\omega}_{X} \subset T_{X} \otimes_{\bold Z} \bold C\, .$$ 
This inclusion defines the Hodge structure of weight $2$ on $T_{X}$. This is a sub-Hodge structure of the Hodge structure of $H^{2}(X, \bold Z)$. We refer to the reader [BPV] for more details.  
Based on the pioneering work of Mukai [Mu], Orlov [Or] proved the following beautiful Theorem: 
\proclaim{Theorem (1.4) [Or Theorem (3.3)]} Let $X_{1}$, $X_{2}$ be K3 surfaces. Then the following three statements are equivalent to one another:

\roster
\item There is a Hodge isometry between $T_{X_{1}}$ and $T_{X_{2}}$: $\varphi_{12} : (T_{X_{1}}, \bold C \omega_{X_{1}}) \rightarrow (T_{X_{2}}, \bold C \omega_{X_{2}})$. 
\item The bounded derived categories of coherent sheaves $D(X_{1})$ and $D(X_{1})$ are equivalent as triangulated category. 

\item $X_{2}$ is isomorphic to (the base space of) a $2$-dimensional compact fine moduli space of some coherent sheaves on $X_{1}$. \qed 
\endroster 
\endproclaim 
Concerning the derived categories, we refer to the reader an excellent book [GM]. Although the equivalence (1) and (3) is implicit in [Or], his proof goes as $\text{(2)} \rightarrow \text{(1)} \rightarrow \text{(3)} \rightarrow \text{(2)}$. See also [BM 1, 2], [Cl], [Th] for further progress and relevant work. 
\flushpar 
\vskip 4pt
(1.5) In their paper [BO], Bondal and Orlov observed that an equivalence class of the derived category $D(V)$ recovers the original manifold $V$ if $V$ is a Fano manifold or a manifold of ample canonical sheaf. Theorem (1.4) says that one can at least recovers the transcendental lattice and the period inside the transcendental lattice if a manifold is a K3 surface. Therefore it is natural to ask: 
\proclaim{Question (1.6)}
\roster 
\item 
how extend does the equivalence class of $D(X)$, or equivalently by (1.4), the Hodge isometry class of $T_{X}$, recover the original K3 surface $X$? 
\item 
at least how extend does the equivalence class of $D(X)$ recover the algebraic part of the lattice $\text{NS}(X)$? \qed
\endroster 
\endproclaim 
In this direction, Mukai [Mu, Proposition (6.2)] already observed that if $\rho(X) \geq 12$, i.e. if $\text{rank}\,T_{X} \leq 10$, then the 
equivalence of $D(X)$ recovers $X$ itself. 
\par
\vskip 4pt 
The aim of this note is to prove the following result: 
\proclaim{Main Theorem (1.7)} Let $N$ be an arbitrarily given natural number. 
Set $I_{N} := \{1, 2, \cdots , N\}$ and $\Delta_{N} := \{(i, i) \vert i \in I\}$. 
Then there are $N$ K3 surfaces $X_{i}$ ($i \in I_{N}$) depending on $N$ such that
\roster 
\item for any $(i, j) \in I_{N}^{2} - \Delta_{N}$, the N\'eron-Severi 
lattices $\text{NS}(X_{i})$ and $\text{NS}(X_{j})$ are not isomorphic, therefore 
$X_{i} \not\simeq X_{j}$ as well, but such that 
\item there is a Hodge isometry 
$\varphi_{ij} : (T_{X_{i}}, \bold C \omega_{X_{i}}) \rightarrow (T_{X_{j}}, \bold C \omega_{X_{j}})$ for any $(i, j) \in I_{N}^{2}$. \qed
\endroster
\endproclaim  
We shall construct such $X_{i}$ of $\rho(X_{i}) = 2$. Note that $2$ is the least possible Picard number for such examples to exist. Indeed, if $\rho(X) = 1$, then $\text{NS}\, (X) = \bold Z l$ and $(l^{2}) = \vert \text{det}T_{X} \vert$ and therefore $\text{NS}\,(X)$ are isomorphic if so are $T_{X}$. 
See (1.10) how extend the Hodge structure of $T_{X}$ recovers $X$ when $\rho(X) = 1$. This is related to the factorization of $\text{deg}\, X$. See also 
Proposition A in the appendix in higher dimensional case.  
\par
\vskip 4pt
Using the global Torelli Theorem of K3 surfaces and Gauss' Theorem of 
quadratic forms, we can reduce the Theorem to a problem of algebraic number theory: find a sequence of real quadratic fields $\bold Q(\sqrt{m})$ having large class number 
$h(m)$ compared with a cardinality $q(m)$ of possible discriminant forms, i.e. $\bold Q/2\bold Z$-valued quadratic forms on the discriminant groups. However, there is no relation between the growth of $h(m)$, $q(m)$ and $m$ even if $m$ is square free. Indeed, there is an open problem whether $h(m) = 1$ for infinitely many $\bold Q(\sqrt{m})$ since Gauss. The growth of $h(m)$ is related to the growth of the norms of the fundamental units and the growth of $q(m)$ is related to the factorization of $m$. 
Theorem (1.2) or the infinite set $\Cal N_{2}$ provides us a way to control the quantities $h(m)$ and $q(m)$ simultaneously when $m = 4n^{2} + 1$ and $n \in \Cal N_{2}$. We shall prove (1.7) in Section 2.  
\par 
As a direct consequence of (1.4) and (1.7), one obtains: 
\proclaim{Corollary (1.8)} For any natural number $N$, 
\roster 
\item there exist at least $N$ K3 surfaces whose derived categories are equivalent but whose Picard groups are non-isomorphic. 
\item there exists a K3 surface $X$ having at least $N$ mutually non-isomorphic 
$2$-dimensional compact fine moduli spaces of stable sheaves on $X$. \qed
\endroster 
\endproclaim 
In the opposite direction, we recall the following fact found by Bridgeland 
and Maciocia (with a bit more constructive proof for the use of the next Proposition (1.10)):
\proclaim{Proposition (1.9) [BK 1, Proposition (5.3)]} Let $X$ be a K3 surface. Let $\Cal S_{T}(X)$ be the set of isomorphism classes of K3 surfaces $Y$ such that there is a Hodge isometry $(T_{Y}, \bold C \omega_{Y}) \simeq (T_{X}, \bold C \omega_{X})$. Then, $\vert \Cal S_{T}(X) \vert < \infty$. 
\qed 
\endproclaim 
(1.9) with (1.4) says that the equivalence class of $D(X)$ recovers K3 
surface $X$ up to finitely many ambiguities as manifold. 
\par 
\vskip 4pt 
It is then interesting to seek explicit values of $\vert \Cal S_{T}(X) \vert$. We shall carry out this in the most generic case, i.e. for K3 surfaces with Picard number 1. Here we write the result in terms of moduli via (1.4):
\proclaim{Proposition (1.10)} Let $X$ be a K3 surface with $\text{NS}(X) = \bold Z l_{X}$. Set $\text{deg}\, X = (l_{X}^{2}) = 2n$. Let $m$ be the number of non-isomorphic $2$-dimemsional compact fine moduli spaces of stable sheaves on $X$. Then $m = \vert ((\bold Z/4n)^{\times})_{2} \vert/4$ if $n \not= 1$ and $m =  \vert (\bold Z/4)^{\times} \vert/2 = 1$ if $n = 1$.  Here $(\bold Z/4n)^{\times}$ is the unit group of the ring $\bold Z/4n$ and $((\bold Z/4n)^{\times})_{2}$ is the two torsion subgroup of $(\bold Z/4n)^{\times}$. More explicitly:
\roster 
\item $m = 1$ if $\text{deg}\,X = 2$; 
\item $m = 1$ if $\text{deg}\,X = 2^{a}$;
\item $m = 2^{k-1}$ if $\text{deg}\,X = 2p_{1}^{e_{1}} \cdots p_{k}^{e_{k}}$; 
\item $m = 2^{k}$ if $\text{deg}\,X = 2^{a}p_{1}^{e_{1}} \cdots p_{k}^{e_{k}}$. 
\endroster 
Here $a$ is a natural number such that $a \geq 2$, $p_{i}$ are mutually different primes such that $p_{i} \geq 3$, and $k$ and $e_{i}$ are natural numbers. 
\qed \endproclaim 
We note that for any given natural number $n$, there is a K3 surface $X$ such that $\rho(X) = 1$ and $\text{deg}\, X = 2n$. 
\head 
{Acknowledgement}
\endhead 
This note has been grown up from a discussion with Professors S. T. Yau, 
S. Hosono, B. Andreas, C. H. Liu at seminar of Professor S. T. Yau. 
Proof owes to a great help of Professor B. H. Gross, the Chair of the Department of Mathematics of Harvard University. The author would like to express his sincere thanks to all of them and to Professors R. P. Thomas 
and E. Kowalski for valuable discussion and comment. The author would like to express his hearty thanks to Professors J. Harris and Y. Kawamata for warm encouragement. This work has been done during his stay at Harvard University under financial support by Harvard University and Japan Monbu-Kagaku Shou (Kaigai-Kenkyu-Kaihatsu-Doukou-Chousa) October 2001 - February 2002. 
\head
{\S 2. Proof of (1.7)}
\endhead 
(2.1) By a {\it lattice} $L := (L, (*, **))$, we mean a free abelian group 
$L$ of finite rank, which we write by $\text{rk}\, L$, equipped with an integral valued non-degenerate symmetric bilinear form $(* ,** ) : L \times L \rightarrow \bold Z$. We write $(x^{2}) = (x, x)$. Two lattices $L_{1} = (L_{1}, (*,**)_{L_{1}})$ and 
$L_{2} = (L_{2}, (*,** )_{L_{2}})$ are said to be {\it isomorphic} or {\it isometric} if there is an isomorphism of abelian groups $f : L_{1} \rightarrow L_{2}$ such that $(f(x), f(y))_{L_{2}} = (x, y)_{L_{1}}$. A lattice $L$ is said to be {\it even} if $(x^{2}) \equiv 0\, \text{mod}\, 2$. We often represent a lattice $L$ by a symmetric integral matrix $S_{L} = ((e_{i}, e_{j}))$ and call $S_{L}$ an associated matrix to $L$. Here 
$\langle e_{i} \rangle_{i=1}^{\text{rk}\, L}$ is integral basis of $L$. An associated matrix is uniquely determined by $L$ up to the action by $\text{GL}(\text{rk}\, L, \bold Z)$: $S \mapsto ^{t}MSM$ on the space of symmetric matrices of degree $\text{rk}\, L$. Note also that two lattices are isomorphic if and only if their associated symmetric matrices lie in the same orbit under the action by $\text{GL}(\text{rk} L, \bold Z)$. We define $\text{det}\, L := \text{det}\, S_{L}$. The signature $(p, q)$ of $L$ is defined to be the signature of an associated matrix $S_{L}$. These are both well-defined and satisfy $q = \text{rk}\, L - p$. We call $L$ {\it hyperbolic} if the signature is $(1, \text{rk} L - 1)$. 
\flushpar
\vskip 4pt
(2.2) We set $L^{*} = \text{Hom}_{\bold Z}(L, \bold Z)$. Using non-degenerate $(*,**)$, one has a natural inclusions $L \subset L^{*} \subset L \otimes_{\bold Z} \bold Q$. Then, $(* ,** )$ is extended to a rational valued symmetric bilinear forms on $L \otimes_{\bold Z} \bold Q$ and on $L^{*}$. Set $A_{L} := L^{*}/L$. We call $A_{L}$ the {\it discriminant group} of $L$. The discriminant group $A_{L}$ is a finite abelian group of order $\vert \text{det} L \vert$. The bilinear form $(*,**)$ of $L$ induces the symmetric bilinear form $B_{A_{L}}(*,**) : A_{L} \times A_{L} \rightarrow \bold Q/ \bold Z$. This is defined by $B_{A_{L}}(\overline{x}, \overline{y}) = (x, y)\, \text{mod}\, \bold Z$, where $x, y \in L^{*}$, $\overline{x} = x\, \text{mod}\, L$, $\overline{y} = y\, \text{mod}\, L$. If $L$ is even, then we have a quadratic form $Q_{A_{L}} : A_{L} \rightarrow \bold Q /2 \bold Z$. This is defined by $Q_{L}(\overline{x}) = (x^{2})\, \text{mod}\, 2\bold Z$. This quadratic form $Q_{A_{L}}$ satisfies 
$$Q_{A_{L}}(n\overline{x}) = n^{2}Q_{A_{L}}(\overline{x})$$ 
$$Q_{A_{L}}(\overline{x} + \overline{y}) 
= Q_{A_{L}}(\overline{x}) + Q_{A_{L}}(\overline{y}) + 2B_{A_{L}}(\overline{x}, \overline{y})\, .$$ 
We call $(A_{L}, Q_{A_{L}})$ a {\it discriminant form} of $L$.
\flushpar 
\vskip 4pt
(2.3) Let $A$ be a finite abelian group. As in [Ni],  
$Q_{A} : A \rightarrow \bold Q/ 2\bold Z$ is called a quadratic form on $A$ if $Q_{A}$ satisfies 
$Q_{A}(na) = n^{2}Q_{A}(a)$ for $a \in A$ and $n \in 
\bold Z$ and if there is a symmetric bilinear map $B(*,**) : A \times A \rightarrow \bold Q / \bold Z$ such that $Q_{A}(a + b) =  Q_{A}(a) + Q_{A}(b) + 2B_{L}(a, b)$. Two elements $a, b \in A$ are said to be orthogonal if $B_{A}(a, b) = 0$. Two pairs $(A, Q_{A})$ and $(A', Q_{A'})$ of finite abelian groups and their quadratic forms are said to be isomorphic if there is a group isomorphism $\varphi : A \rightarrow A'$ such that $Q_{A}(a) = Q_{A'}(\varphi(a))$ for all $a \in A$. We define the orthogonal sum $(A_{1} \oplus A_{2}, Q_{A_{1}} \oplus Q_{A_{2}})$ of $(A_{1}, Q_{A_{1}})$ and $(A_{2}, Q_{A_{2}})$ in a natural manner: That is,  $A_{1} \oplus A_{2}$ is the direct sum as abelian group and $Q_{A_{1}} \oplus Q_{A_{2}}(a_{1} \oplus a_{2}) := Q_{A_{1}}(a_{1}) + Q_{A_{2}}(a_{2})$.  
\par 
\vskip 4pt
Let $L$ be an even lattice. One can decompose $A_{L}$ into the direct sum of the Sylow subgroups of $A_{L}$: 
$$A_{L} = \oplus_{i = 1}^{k} (A_{L})_{p_{i}}\, .$$
Here $(A_{L})_{p_{i}}$ is the Sylow $p_{i}$-subgroup of $A_{L}$. 
\proclaim{Lemma (2.4) [Ni, Page 108-109]} 
\roster 
\item If $i \not= j$, then $B_{L}(\overline{x}_{i}, \overline{x}_{j}) = 0$ for $\overline{x}_{i} \in (A_{L})_{p_{i}}$ and $\overline{x}_{j} \in (A_{L})_{p_{j}}$. 
\item Let $x \in L^{*}$ and write $\overline{x} = \sum_{i=1}^{k}\overline{x}_{i}$, where $\overline{x}_{i} \in (A_{L})_{p_{i}}$. Then, $Q_{L}(\overline{x}) = \sum_{i=1}^{k}Q_{L}(\overline{x}_{i})$. In other words, we have the orthogonal decomposition 
$$(A_{L}, Q_{L}) = \oplus_{i=1}^{k} ((A_{L})_{p_{i}}, Q_{L} \vert (A_{L})_{p_{i}}).$$ 
\endroster
\endproclaim  
\demo{Proof} It is clear that (1) implies (2). Let us show (1). Take $\overline{x}_{i} \in (A_{L})_{p_{i}}$ and $\overline{x}_{j} \in (A_{L})_{p_{j}}$ and set $\text{ord}\,(\overline{x}_{i}) = p_{i}^{n}$. Since $p_{i}^{n}$ and $p_{j}$ are coprime, there is $\overline{y}_{j} \in (A_{L})_{p_{j}}$ such that $\overline{x}_{j} = p_{i}^{n}\overline{y}_{j}$. Then one calculates 
$$B_{L}(\overline{x}_{i},\overline{x}_{j}) 
= B_{L}(\overline{x}_{i}, p_{i}^{n}\overline{y}_{j}) = 
B_{L}(p_{i}^{n}\overline{x}_{i}, \overline{y}_{j}) 
= B_{L}(0, \overline{y}_{j}) = 0\, .\, \qed$$ \enddemo
\proclaim{Lemma (2.5)} Let $p \not= q$ be odd primes. Set 
$\Cal L_{p, q}$ to be the isomorphism classes of even hyperbolic lattices 
$L$ such that $\vert \text{det}\, L \vert = pq$. Set 
$\Cal Q_{p, q}$ to be the set of the isomorphism classes of $(A_{L}, Q_{A_{L}})$ such that $L \in \Cal L_{p, q}$. Then there is a positive integer $B$ being independent of $p, q$ such that $\vert \Cal Q_{p, q} \vert \leq B$ for all $p \not= q$. 
\endproclaim 
\demo{Proof} Let $L \in \Cal L_{p, q}$. Then the decomposition of $A_{L}$ into the Sylow subgroups is $\bold Z/p \oplus \bold Z/q$. By (2.4), one has 
$Q_{L} = Q_{p} \oplus 
Q_{q}$, where $Q_{p}$ and $Q_{q}$ are quadratic forms on $\bold Z/p$ and $\bold Z/q$ respectively. It suffices to estimate the cardinality of the isomorphism classes $(\bold Z/p, Q_{p})$. Note that $Q_{p}$ is determined by the value $Q_{p}(1)$. 
\par 
If $p = 1$, then there is just one quadratic form on $\bold Z/ 1 = \{0\}$. So, we may assume $p \geq 3$. Choose and fix an element $p_{-} \in \bold N$ such that 
$(\frac{p_{-}}{p}) = -1$. Here $(\frac{*}{**})$ is the Legendre symbol of the quadratic residue. By $Q_{p}(p) = Q_{p}(0) = 0$, one has $Q_{p}(1) = 2a/p$, where $a$ is some integer. Note that the maps $n : \bold Z/p \rightarrow \bold Z/p$ defined by $1 \mapsto n \cdot 1$ are all isomorphism as abelian group provided that $n$ is coprime to $p$ and satisfy $Q_{p}(n) = n^{2}Q_{p}(1)$. 
Thus, $(\bold Z/p, Q_{p})$ is isomorphic to either $(\bold Z/p, q_{0})$, 
$(\bold Z/p, q_{+})$ or $(\bold Z/p, q_{-})$, where $q_{*}$ are defined by $q_{0}(1) = 0$, $q_{+}(1) = 2/p$, $q_{-}(1) = 2p_{-}/p$. Therefore one can 
take $B = 3 \cdot 3 = 9$. \qed \enddemo
Let $m > 1$ be a square free natural number such that $m \equiv 1 (4)$. Then $\bold Q(\sqrt{m})$ is a real quadratic field.  We write $\bold Q(\sqrt{m})$ by $K(m)$. Denote by $O(m)$, $h(m)$, $\epsilon(m) > 1$, $D(m)$ the ring of integers, the class number, the fundamental unit, the discriminant of $K(m)$ respectively. $O(m)$ is by definition, the normalization of $\bold Z$ in $\bold K(m)$ and $O(m) = \bold Z[(1 + \sqrt{m})/2]$. One has $D(m) = m$. The unit group $O(m)^{\times}$ of $O(m)$ satisfies $O(m)^{\times} = \langle \epsilon(m) \rangle \times \langle -1 \rangle \simeq \bold Z \times \bold Z/2$. Among the four possible free generators 
$\pm \epsilon(m)^{\pm 1}$ of $O(m)^{\times}$, the fundamental unit is the one which is greater than $1$. The class number $h(m)$ is defined to be the cardinality of the ideal class group $C(m)$ of $O(m)$. We need the following two Theorems (2.7) and (2.8).
The first one is attributed to Gauss and the second one is (a special case of) the Theorem of Siegel-Brauer. For the statement, we recall that an integral binary quadratic form $f(x, y) = ax^{2} + bxy + cy^{2}$ is called {\it primitive} if $(a, b, c) = 1$. The integer $d(f) := b^{2} - 4ac$ is called the discriminant of $f$. By $S_{f}$, we denote the associated symmetric (rational) matrix, i.e. 
$$S_{f} = \pmatrix a & \frac{b}{2}\\
            \frac{b}{2} & c\\ \endpmatrix .$$ 
Note that $\text{det}\,S_{f} = -d(f)/4$. We call two integral binary quadratic forms $f$ and $f'$ are {\it properly equivalent} if $S_{f}$ and $S_{f'}$ lie in the same orbit under the natural action by $\text{SL}(2, \bold Z)$, i.e. if there is $M \in \text{SL}(2, \bold Z)$ such that $^{t}MS_{f}M = S_{f'}$. 
\proclaim{Theorem (2.7) (e.g. [Nr, Theorem (8.6)])} There is a natural one-to-one correspondence between the class group $C(m)$ and the set of the properly equivalent classes of primitive integral binary forms of discriminant $m$. 
In particular, $h(m)$ is the cardinality of the set of the properly equivalent classes of primitive integral binary forms of discriminant $m$. 
\qed \endproclaim   
\proclaim{Theorem (2.8) (e.g. [Nr, Theorem (8.5)])} 
$$\lim_{m \rightarrow \infty}\frac{\log (h(m) \cdot \log \epsilon(m))}{\log D(m)} = \frac{1}{2}\, .\, \qed$$
\endproclaim 
The following Proposition is crucial for the main Theorem:
\proclaim{Proposition (2.9)} Let $N$ be an arbitrarily given natural number. 
Set $I_{N} := \{1, 2, \cdots , N\}$ and $\Delta_{N} := \{(i, i) \vert i \in I_{N}\} \subset I_{N}^{2}$. Then, there exist 
$N$ even hyperbolic lattices of rank $2$, say, $S_{i}$ ($i \in I_{N}$) such that 
\roster 
\item if $(i, j) \in I_{N}^{2} - \Delta_{N}$, then $S_{i}$ and $S_{j}$ are not isomorphic, but
\item for all $(i, j) \in I_{N}^{2}$, the discriminant forms $(A_{S_{i}}, Q_{A_{S_{i}}})$ and $(A_{S_{j}}, Q_{A_{S_{j}}})$ are isomorphic. 
\endroster 
\endproclaim 
\demo{Proof} Let $\Cal N_{2}$ be the set in the Introduction. Let $n \in \Cal N_{2}$ and write $d(n) := 4n^{2} + 1 = pq$. 
Note that $d(n) \equiv 1\, \text{mod}\, 4$. Set $K_{n} := K(d(n))$, 
$O_{n} := O(d(n))$, $O_{n}^{\times} = O(d(n))^{\times}$, $\epsilon_{n} := 
\epsilon(d(n))$, $h_{n} := h(d(n))$, $D_{n} := D(d(n)) = 4n^{2} + 1 = pq$ under the notation in (2.6). Recall that $\Cal N_{2}$ is an infinite set. 
\enddemo
\proclaim{Lemma (2.10)} 
For any given $\epsilon > 0$, there is a natural number $M$ such that 
$$\frac{\log (\log \epsilon_{n})}{\log D_{n}} < \epsilon$$
for all natural numbers $n$ such that $n > M$ and $n \in \Cal N_{2}$. 
\endproclaim 
\demo{Proof of (2.10)} Consider $2n + \sqrt{pq} \in O_{n}$. We may assume that $n \geq 3$. We have $\log \epsilon_{n} > 0$ by the definition. Note that
$$(2n + \sqrt{pq})(2n - \sqrt{pq}) = 4n^{2} - pq = 4n^{2} -(4n^{2} +1) = -1\, .$$ 
From this identity, we obtain 
$$\frac{1}{2n + \sqrt{pq}} = - (2n - \sqrt{pq}) \in O_{n}\, .$$
In particular, $2n + \sqrt{pq} \in O_{n}^{\times}$. Combining this with $2n + \sqrt{pq} > 1$, we find a natural number $l$ such that $(2n + \sqrt{pq})^{2} = \epsilon_{n}^{l}$. Using $4n \geq \sqrt{4n^{2} + 1} = \sqrt{pq}$, $n \geq 3$ and $\epsilon_{n} > 1$, we calculate that 
$$4 \log 2n \geq 2 \log 2n + 2 \log 3 = 2 \log 6n \geq 2 
\log (2n + \sqrt{pq})$$ 
$$= \log (2n + \sqrt{pq})^{2} = \log \epsilon_{n}^{l} = l \log \epsilon_{n} \geq \log 
\epsilon_{n}\, .$$ 
We have also 
$$D_{n} = 4n^{2} + 1 \geq 2n\,.$$
Combining these two inequalities, one has 
$$\frac{\log (\log \epsilon_{n})}{\log D_{n}} \leq \frac{\log (4 \log 2n)}
{\log 2n} = \frac{\log(\log 2n) + \log 4}
{\log 2n}\, .$$
Since $\lim_{x \rightarrow \infty} \log x = \infty$ and $\lim_{x \rightarrow \infty} (\log x)/x = 0$, this gives the result. \qed 
\enddemo 
Combining (2.8) with (2.10) applied for $\epsilon = 1/8$, we obtain 
\proclaim{Corollary (2.11)} There is a natural number $M$ such that 
$$\frac{\log h_{n}}{\log D_{n}} > \frac{1}{4}$$
for all natural numbers $n$ such that $n > M$ and $n \in \Cal N_{2}$. 
In particular, 
$$\lim_{n \rightarrow \infty, n \in \Cal N_{2}} h_{n} = \infty\, .\, \qed$$ 
\endproclaim
Let us return back to the proof of (2.9). Let $B$ be the constant in (2.5). 
Take an arbitrary natural number $N$. By (2.11), there is $n \in \Cal N_{2}$ such that $h_{n} > 2BN$. We write $4n^{2} + 1 = pq$. 
Then by (2.7), there are more than $2BN$ binary forms of discriminant $pq$ which are not properly isomorphic to one another. We write them by 
$a_{i}x^{2} + b_{i}xy + c_{i}y^{2}$, $i \in \{1, 2, \cdots , 2BN \}$. These binary forms satisfy $b_{i}^{2} - 4a_{i}c_{i} = pq$. Then, the associated symmetric matrices 
multiplied by $2$, i.e. 
$$M_{1} = \pmatrix 2a_{1} & b_{1}\\ b_{1} & 2c_{1}\\ \endpmatrix\, , \cdots , 
M_{2BN} = \pmatrix 2a_{2BN} & b_{2BN}\\ b_{2BN} & 2c_{2BN}\\ \endpmatrix$$ 
are all even, hyperbolic, and of determinant $-pq$ and lie in mutually different orbits under the natural action by $\text{SL}(2, \bold Z)$.  Since $\text{SL}(2, \bold Z)$ 
is a normal subgroup of $\text{GL}(2, \bold Z)$ of index $2$, among these $2BN$ matrices $M_{i}$, one can find $BN$ matrices lying in mutually different orbits under the action by $\text{GL}(2, \bold Z)$. By renumbering, we may assume they are $M_{i}$ $(i \in \{1, \cdots , BN\})$. Then, these $BN$ matrices $M_{i}$ ($i \in \{1, \cdots , BN\}$) define $BN$ mutually non-isomorphic even hyperbolic lattices $S_{i}$ of determinant $-pq <0$ and of rank $2$. On the other hand, by (2.5), there are at most $B$ discriminant forms $(A_{S_{i}}, Q_{A_{S_{i}}})$ up to isomorphism when $i$ runs through $\{1, \cdots , BN\}$. Therefore, among these $S_{i}$ ($i \in \{1, \cdots , BN\}$), there are at least $BN/B = N$ lattices, say, $S_{i}$ ($i \in \{1, \cdots , N\}$), which have mutually isomorphic discriminant forms. Now we are done for (2.9). \qed
\par
\vskip 4pt
We shall construct K3 surfaces in the main Theorem (1.7) by using the lattices  $S_{i}$ in (2.9). For this, we need the following Theorem due to Nikulin. For the statement, we recall that an embedding as lattice $\Phi : M \rightarrow L$ is said to be {\it primitive} if the abelian group $L/\Phi(M)$ is free, in other words, free basis of $\Phi(M)$ can be extended to free basis of $L$, or passing to the dual, the natural homomorphism $L^{*} \rightarrow \Phi(M)^{*}$ is surjective.  A sublattice $M$ of $L$ is called primitive if the inclusion map is primitive.
\proclaim{Theorem (2.11) [Ni, Theorem (1.14.4), Corollary (1.13.3)]} Let $M$ be an even lattice of signature $(m_{+}, m_{-})$ (therefore $\text{rk} \,M = m_{+} + m_{-}$). Let $L$ be an even unimodular lattice of signature $(l_{+}, l_{-})$.
\roster 
\item Let $\Phi : M \rightarrow L$ be a primitive embedding and $K$ the orthogonal lattice of $\Phi(M)$ in $L$, i.e. $K = \{x \in L \vert (x, y) = 0 \,\,\text{for all}\,\, y \in \Phi(M)\}$. (Note that $K$ is primitive in $L$.) Then, under isomorphism 
$$M^{*}/M \simeq \Phi(M)^{*}/\Phi(M) \simeq L/(\Phi(M) \oplus K) \simeq K^{*}/K$$ given by the natural surjective homomorphism $L^{*} = L \rightarrow \Phi(M)^{*}$ and  $L^{*} = L \rightarrow K^{*}$, one has 
$$(A_{M}, Q_{A_{M}}) \simeq (A_{K}, -Q_{A_{K}})\, .$$ 
\item Assume that $l_{+} - t_{+} > 0$, $l_{-} - t_{-} > 0$ and $\text{rk}\, L - \text{rk}\, M \geq 2 + l(A_{M})$, where $l(A_{M})$ is the minimal number of generators of the finite abelian group $A_{M} = M^{*}/M$. Then, $M$ can be primitively embedded into $L$. Moreover, a primitive embedding $M \rightarrow L$ is unique in the sense that if $f_{i} : M \rightarrow L$ ($i = 1, 2$) are two primitive embeddings, then there is an isometry $\varphi : L \rightarrow L$ such that $f_{2} = \varphi \circ f_{1}$. 
\item Let $\tilde{M}$ be an even lattice. Assume that the signature of 
$\tilde{M}$ is $(m_{+}, m_{-})$, i.e. the same as the signature of $M$ 
and that $(A_{\tilde{M}}, Q_{A_{\tilde{M}}})$ is isomorphic to $(A_{M}, Q_{A_{M}})$. Then $\tilde{M}$ is isomorphic to $M$ provided that $\text{rk}\, M \geq 2 + l(A_{M})$, $m_{+} > 0$ and $m_{-} > 0$.\, \qed 
\endroster
\endproclaim  
Set $\Lambda := U^{\oplus 3} \oplus E_{8}(-1)^{\oplus 2}$. Here $U$ is the lattice given by $\pmatrix 0 & 1\\ 1 & 0\\ \endpmatrix$ and $E_{8}(-1)$ is 
the negative lattice given by the Dynkin diagram of type $E_{8}$. This $\Lambda$ is an even unimodular lattice of signature $(3, 19)$ and is called the {\it K3 lattice}. This is a unique even unimodular lattice of signature $(3, 19)$ up to isomorphism. For any K3 surface $X$, there is an isometry $\tau : H^{2}(X, \bold Z) \rightarrow \Lambda$ (See for instance [BPV]). We call an isometry $\tau$ from $H^{2}(X, \bold Z)$ to $\Lambda$ a marking of $X$. We call such a pair $(X, \tau)$ a {\it marked K3 surface}. 
\proclaim{Lemma (2.12)} Let $S_{i}$ ($1 \leq i \leq N$) be the lattices found in (2.9). 
Then 
\roster
\item $S_{i}$ are primitively embedded into $\Lambda$, say, $\Phi_{i} : S_{i} \rightarrow \Lambda$. 
\item Denote the orthogonal lattices of $\Phi_{i}(S_{i})$ in $\Lambda$ by $T_{i}$. Then $T_{i}$ are isomorphic to one another. 
\endroster
\endproclaim 
\demo{Proof} By the elementary divisor theory, one has $l(A_{S}) \leq \text{rk}\, S$. Since $\text{rk}\, S_{i} = 2$, (2.11)(2) implies (1).  Let us check (2). Since $(A_{S_{i}}, Q_{A_{S_{i}}}) \simeq (A_{S_{j}}, Q_{ A_{S_{j}}})$ by the construction, we have $(A_{T_{i}}, Q_{A_{T_{i}}}) \simeq 
(A_{T_{j}}, Q_{A_{T_{j}}})$ by (2.11)(1). Since $S_{i}$ is of signature $(1,1)$ and 
$\Lambda$ is of signature $(3, 19)$, $T_{i}$ are of signature $(2, 18)$.  By (2.11)(1), we have also $l(A_{T_{i}}) = l(A_{S_{i}}) \leq 2$. Now the assertion (2) follows from (2.11)(3). \qed
\enddemo
Let $T_{i}$ ($i \in \{1, \cdots , N\}$) be the sub-lattices of $\Lambda$ found in (2.12). Take an even lattice $T$ of signature $(2, 18)$ which is isomorphic to all these $T_{i}$. We denote by $\varphi_{i} : T \rightarrow \Lambda$ a primitive embedding of $T$ such that $\varphi_{i}(T) = T_{i}$. 
\proclaim{Lemma (2.13)} There is a positive definite $2$-dimensional subspace $P$ of $T \otimes_{\bold Z} \bold R$ such that $P^{\perp} \cap T = \{0\}$, where $P^{\perp}$ is the orthogonal space of $P$ in $T \otimes_{\bold Z} \bold R$. 
\endproclaim 
\demo{Proof} Since the signature of $T$ is $(2, 18)$, there is a positive definite $2$-dimensional subspace $P_{0}$ of $T \otimes_{\bold Z} \bold R$. Note that positive definiteness is an open condition in the real Grassman manifold $\text{Gr}(2, T \otimes_{\bold Z} \bold R)$ in the classical topology. Moreover, in $\text{Gr}(2, T \otimes_{\bold Z} \bold R)$, the locus $P^{\perp} \cap T \not= \{0\}$ is countable union of the proper Zariski closed subsets $(P, t) = 0$ ($t \in T - \{0\}$) in $\text{Gr}(2, T \otimes_{\bold Z} \bold R)$. Here the properness is because $T$ is non-degenerate. Therefore, there is a desired $P$ near $P_{0}$ in $\text{Gr}(2, T \otimes_{\bold Z} \bold R)$. 
\qed
\enddemo 
(2.14) Let us complete the proof of (1.7). Let $P$ be the space found in (2.13). Let $\langle \eta_{1}, \eta_{2} \rangle$ be the orthonomal basis of $P$. Set $\omega := \eta_{1} + \sqrt{-1}\eta_{2}$ in $P \otimes_{\bold R} \bold C$. Then, by $(\eta_{1}, \eta_{1}) = (\eta_{2}, \eta_{2}) = 1$ and $(\eta_{1}, \eta_{2}) = 0$, one has $(\omega, \omega) = 0$ and 
$(\omega, \overline{\omega}) = 2 > 0$. Here $\overline{\omega}$ is the complex conjugate of $\omega$ with respect to the real structure $P$. Therefore if we set $\omega_{i} = \varphi_{i}(\omega) \in \Lambda \otimes_{\bold Z} \bold C$, then $\omega_{i} \in T_{i} \otimes_{\bold Z} \bold C$, $(\omega_{i}, \omega_{i}) = 0$ and $(\omega_{i}, \overline{\omega_{i}}) > 0$ as well. Moreover, by (2.13) and by $\Phi_{i}(S_{i}) = T_{i}^{\perp}$ in $\Lambda$, the following equalities hold in $\Lambda \otimes_{\bold Z} \bold C$:
$$\omega_{i}^{\perp} \cap \Lambda 
= \langle \varphi_{i}(\eta_{1}), \varphi_{i}(\eta_{2}) \rangle ^{\perp} \cap \Lambda$$ 
$$= \varphi_{i}(P)^{\perp} \cap \Lambda = T_{i}^{\perp} \cap \Lambda = \Phi_{i}(S_{i}).$$ 
Therefore $(\Lambda, \bold C \omega_{i})$ ($i \in \{1, \cdots , N\}$) define the weight two Hodge structures on $\Lambda$ such that $\omega_{i}^{\perp} = \Phi_{i}(S_{i})$. Then, by the surjectivity of the period mapping of K3 surfaces and by the Lefschetz $(1, 1)$-Theorem, there are marked K3 surfaces $(X_{i}, \tau_{i})$ such that $\tau(\bold C \omega_{X_{i}}) = \bold C \omega_{i}$ and $\tau_{i}(\text{NS}(X_{i})) = \Phi_{i}(S_{i})$. 
This implies $\tau_{i}(T_{X_{i}}) = T_{i}$ as well. 
Since $\Phi_{i}(S_{i})$ are hyperbolic, so are $NS(X_{i})$. Therefore $X_{i}$ are projective. Moreover, by the construction, the following homomorphism gives a Hodge isometry between $(T_{X_{i}}, \bold C \omega_{X_{i}})$ and $(T_{X_{j}}, \bold C \omega_{X_{i}})$: 
$$T_{X_{i}} \overset{\tau_{i}} \to \longrightarrow 
\varphi_{i}(T_{i}) \overset{\varphi_{i}^{-1}} \to \longrightarrow 
T \overset{\varphi_{j}} \to \longrightarrow 
\varphi_{j}(T_{j}) \overset{\tau_{j}^{-1}} \to \longrightarrow 
T_{X_{j}}\, .$$
These K3 surfaces $X_{i}$($i \in \{1, \cdots , N\}$) satisfy all the requirement in (1.7). \qed  
\head 
{\S 3. Proof of (1.9)}
\endhead 
Let $(T, \bold C \omega)$ be a lattice with weight two Hodge structure isomorphic to the Hodge structure $(T_{X}, \bold C \omega_{X})$. (For this notation, we note that the Hodge structure on the lattice $T_{X}$ is determined by the inclusion $\bold C \omega_{X} \subset T_{X} \otimes_{\bold Z} \bold C$.) 
Let $Y \in \Cal S_{T}(X)$. Then $\text{rk}\, NS(Y) = 22 - \text{rk}\, T$ and $\vert \text{det}\, \text{NS}\,(Y) \vert = \vert \text{det}\,T_{Y} \vert = \vert \text{det}\, T_{Y} \vert$. Therefore, the finiteness of reduction of non-degenerate integral quadratic forms with bounded determinant and bounded rank (see for instance [Cs, Page 128, Theorem 1.1]), one has 
$$\vert \{\text{NS}(Y) \vert Y \in \Cal S_{T}(X)\}/\text{isom} \vert < \infty\, .$$ 
Let $S_{i}$ ($i \in \{1, \cdots , N\}$) be the complete representative of the set above. 
Then one has $\Cal S_{T}(X) = \cup_{i=1}^{N} \Cal S_{i}(X)$ (disjoint union), where 
$$\Cal S_{i}(X) := \{Y \in \Cal S_{T}(X) \vert \text{NS}\,(Y) \simeq S_{i}\}/\text{isom}\, .$$ 
Therefore, it suffices to show that $\vert \Cal S_{i}(X) \vert < \infty$ for each $i$. Choose $i$. For simplicity of notation, we write $S$, $\Cal S(X)$ for $S_{i}$, $\Cal S_{i}(X)$ and so on from now. 
\par 
\vskip 4pt 
A $\bold Z$-module $L$ such that $S \oplus T \subset L \subset S^{*} \oplus T^{*}$ is called an {\it over lattice} of $S \oplus T$. 
Let us consider all the even unimodular over lattices $L$ of $S \oplus T$ such that $S$ and $T$ are both primitive in $L$. Such an $L$ is an even unimodular lattice of signature $(3, 19)$. Therefore $L$ is isomorphic to the K3 lattice. Since $(S^{*} \oplus T^{*})/(S \oplus T) = S^{*}/S \oplus T^{*}/T$ is a finite group, there are only finitely many such $L$ as a subset of $S^{*} \oplus T^{*}$. We write all of them by $L_{j}$ ($j \in \{1, \cdots , M\}$). 
\par 
\vskip 4pt
Let $Y \in \Cal S(X)$. Then, we have an isometry $f : \text{NS}(Y) \rightarrow S$ and a Hodge isometry $g : (T_{Y}, \bold C \omega_{Y}) \rightarrow (T, \bold C \omega)$. 
This induces an isomorphism $(f \oplus g) : S_{Y}^{*} \oplus T_{Y}^{*} \rightarrow S^{*} \oplus T^{*}$. Then there is $j \in \{1, \cdots , M\}$ such that $(f \oplus g)(H^{2}(Y, \bold Z)) = L_{j}$. Write $\tau_{Y} := (f \oplus g) \vert H^{2}(Y, \bold Z)$. Conversely, by the surjectivity of the period mapping of K3 surfaces, for each $j \in \{1, \cdots , M\}$, there is a marked K3 surface $(X_{j}, \tau_{X_{j}})$ such that $\tau_{X_{j}} : H^{2}(X_{j}, \bold Z) \rightarrow L_{j}$ is a Hodge isometry in the sense that $\tau_{j}(\omega_{X_{j}}) = \omega$. Since $T$ and $S$ are both primitive in $L_{j}$, 
we have $\tau_{X_{j}}(T_{X_{j}}) = T$ and $\tau_{X_{j}}(\text{NS}(X_{j})) 
= S$ in this case. Moreover, if we put $f = \tau_{X_{j}}\vert \text{NS}(X_{j})$ and $g = \tau_{X_{j}} \vert T_{X_{j}}$, then $\tau_{X_{j}}$ is recovered from $f$ and $g$ by the process explained above. Let us choose for each $j \in \{1, \cdots , M\}$ a marked K3 surface $(X_{j}, \tau_{X_{j}})$ as above. 
Let $Z \in \Cal S(X)$. Then $\tau_{Z}(H^{2}(Z, \bold C)) = L_{j}$ for some $j \in \{1, \cdots , M\}$. By construction, we have 
$\tau_{Z}(H^{2}(Z, \bold Z)) = L_{j} = \tau_{X_{j}}(H^{2}(Y_{j}, \bold Z))$ and $\tau(\bold C \omega_{Z}) = \bold C \omega = \tau_{X_{j}}(\bold C \omega_{X_{j}})$. Thus, $\tau_{Z}^{-1} \circ \tau_{X_{j}} : 
H^{2}(X_{j}, \bold Z) \rightarrow H^{2}(Z, \bold Z)$ is a Hodge isometry. By the global Torelli Theorem for K3 surface, we have then $Z \simeq X_{j}$. Therefore $\Cal S(X)$ consists of at most $M$ elements. \qed 
\head 
{\S 4. Proof of (1.10)}
\endhead 
Idea of proof of (1.10) is similar to that of (1.9), but we need a bit more precise argument in order to obtain the exact number $\vert \Cal S_{T}(X) \vert$.
\par
\vskip 4pt
Let $S = \langle l \rangle$ be a lattice of rank $1$ such that $(l^{2}) = 2n$. Let $X$ be the same as in (1.10). As before, we choose an abstract lattice with weight two Hodge structure $(T, \bold C \omega)$ isomorphic to $(T_{X}, \bold C \omega_{X})$. If $Y \in \Cal S_{T}(X)$, then $\text{NS}(Y) = \bold Z l_{Y}$ and $(l_{Y})^{2} = \vert \text{det} T_{Y} \vert = \vert \text{det} T_{X} \vert = (l_{X}^{2}) = 2n$. We take $l_{Y}$ the ample class. So, $l_{Y}$ is uniquely determined by $Y$. We have an isometry $f_{Y} : \text{NS}(Y) \simeq S; l_{Y} \mapsto l$ and a Hodge isometry $g_{Y} : (T_{Y}, \bold C \omega_{Y}) \rightarrow (T, \bold C \omega)$.  We choose and fix $f_{Y}$ and $g_{Y}$ for each $Y \in \Cal S_{T}(X)$.  
\proclaim{Lemma (4.1)} Hodge isometry from $(T_{Y}, \bold C \omega_{Y})$ to $(T, \bold C \omega)$ is either $g_{Y}$ or $-g_{Y}$. Conversely both are Hodge isometries. 
\endproclaim 
\demo{Proof} The last statement is clear. Let us show the first assertion. Let $g : (T_{Y}, \bold C \omega_{Y}) \rightarrow (T, \bold C \omega)$ be a Hodge isometry. Then $g_{Y}^{-1} \circ g$ is a Hodge isometry of $(T_{Y}, \bold C \omega_{Y})$. So, it suffices to show that if $h$ is a Hodge isometry of $(T_{Y}, \bold C \omega_{Y})$, 
then $h = id$ or $-id$. 
\par 
Since $T \otimes_{\bold Z} \bold R = P \oplus N$, where 
$P = (\bold C \omega_{Y} \oplus \bold C \overline{\omega_{Y}}) \cap (T\otimes_{\bold Z} \bold R)$ and $N = H^{1,1}(Y, \bold C) \cap (T\otimes_{\bold Z} \bold R)$. Here $P$ is positive definite. Since $Y$ is projective, $N$ is negative definite by the Hodge index Theorem. 
Since $h$ is a Hodge isometry, $h \in O(P) \times O(N)$. Here $O(*)$ is the orthogonal group. Since $P$ and $N$ are both definite, $h$ is diagonalizable and the eigenvalues of $h$ (in $\bold C$) are all absolute value $1$. On the other hand, since $h$ is defined over $\bold Z$, the eigenvalues of $h$ are all algebraic integers. Therefore, the eigenvalues of $h$ are all root of unity by the Theorem of Kronecker. In particular, there is a natural number $I$ such that $h(\omega_{Y}) = \zeta_{I}\omega_{Y}$, where $\zeta_{I}$ is a primitive $I$-th root of unity. Then $\text{ord}(h) = I$. Otherwise, we have $\text{ord}(h) = kI$ for some integer $k \geq 2$. However, since $h$ is defined over $\bold Z$, the space $T'= \{x \in T_{Y} \vert h^{I}(x) = x \}$ would be a primitive sub-module of $H^{2}(Y, \bold Z)$, $T' \not= T$ and $\omega_{Y} \in T' \otimes_{\bold Z} \otimes C$. However, this contradicts the direct consequence of the Lefschetz $(1,1)$-Theorem: $T_{Y}$ is the minimal primitive sub-module of $H^{2}(Y, \bold Z)$ such that $\omega_{Y} \in T_{Y} \otimes_{\bold Z} \bold C$.     
\par 
On the other hand, since $\text{dim}\, N = \text{rk}\, T_{Y} -2 = 19$ is odd, $h \vert N$ and therefore $h$ has $1$ or $-1$ as its eigenvalue. Since $h$ is defined over $T_{Y}$ and $\pm 1$ is rational, there then exists an element $a \in T_{Y} - \{0\}$ such that either $h(a) = a$ or $h(a) = -a$. Therefore we have two cases: 
\roster 
\item $h$ has eigenvalue $1$; 
\item $h$ does not have eigenvalue $1$ but has eigen value $-1$. 
\endroster 
We shall show that in the first case $h = id$ and in the second case $h = -id$. \par 
Let us consider the Case 1. Take $a \in T_{Y} -\{0\}$ such that $h(a) = a$. 
Then:
$$(\zeta_{I}\omega_{Y}, a) = (h(\omega_{Y}), h(a)) = (\omega_{Y}, a)\, .$$ 
If $I \not= 1$, this would imply $(\omega_{Y}, a) = 0$. Then $a \in \text{NS}(Y) \cap T_{Y} = \{0\}$, a contradiction. Therefore $I = 1$ and 
$h = id$. 
\par 
Let us consider the Case 2. There is $a \in T_{Y} -\{0\}$ such that $h(a) = -a$.  Then the formula 
$$(\zeta_{I}\omega_{Y}, -a) = (h(\omega_{Y}), h(a)) = (\omega_{Y}, a)$$
implies $I = 2$. Indeed, otherwise, we would have $(\omega_{Y}, a) = 0$ and 
get the same contradiction as Case 1. From this and the case assumption, the eigenvalues of $h$ are all $-1$. Since $h$ is diagonalizable, this implies 
$h = -id$. \qed \enddemo
As in Section 3, we take all the even unimodular over lattices $L$ 
of $S \oplus T$ such that $S$ and $T$ are both primitive in $L$. Write all of them by $L_{j}$ ($j \in \{1, \cdots , M\}$). 
As in the proof of (1.9), if $Y \in \Cal S_{T}(X)$, then under the natural extension $(f_{Y} \oplus g_{Y}) : \text{NS}(Y)^{*} \oplus T_{Y}^{*} \rightarrow S^{*} \oplus T^{*}$, there is $j$ such that 
$(f_{Y} \oplus g_{Y})(H^{2}(Y, \bold Z)) = L_{j}$. We write $\tau_{Y} := (f_{Y} \oplus g_{Y})\vert H^{2}(Y, \bold Z)$ as before. We already observed the following fact in Section 3:
\proclaim{Lemma (4.2)} 
\roster 
\item 
For each $j \in \{1, \cdots , M\}$, there is a K3 surface 
$X_{j} \in \Cal S_{T}(X)$ such that $\tau_{X_{j}}(H^{2}(X_{j}, \bold Z)) = L_{j}$. 
We fix such an $X_{j}$ for each $j$ in what follows.
\item 
For each $Y \in \Cal S_{T}(X)$, there is $j \in \{1, \cdots , M\}$ such that 
$Y \simeq X_{j}$. \qed
\endroster 
\endproclaim 
So, one can find complete representatives of $\Cal S_{T}(X)$ in $\{X_{j}\}_{j=1}^{M}$. 
\par 
\vskip 4pt
However, by (4.1), there are exactly two choices of $g_{X_{j}}$ for each $X_{j} \in \Cal S_{T}(X)$: $g_{X_{j}}$ and $-g_{X_{j}}$. Therefore in order to find out the complete representatives in $\{X_{j}\}_{j=1}^{M}$, we also need to 
seek how $L_{j}$ changes when we replace $g_{X_{j}}$ by $-g_{X_{j}}$. 
\par 
\vskip 4pt 
Let $L$ be any one of $L_{j}$. Then $T^{*}/T \simeq L/(T \oplus S) \simeq S^{*}/S \simeq \bold Z/2n$ by (2.11)(1). Therefore, $S^{*}/S = \langle \frac{l}{2n} \rangle$ 
and there is $t \in T$ such that $T^{*}/T =  \langle \frac{t}{2n} \rangle$. 
This $t$ is independent of $L$. Note that $L/(T \oplus S) = \langle \frac{al + bt}{2n} \rangle$ for some $a, b \in \bold Z$. This is because $\bold Z/2n \simeq L/(T \oplus S) \subset S^{*}/S \oplus T^{*}/T = (S^{*} \oplus T^{*})/(S \oplus T)$. Note also that under the natural isomorphism $L/(T \oplus S) \simeq S^{*}/S$ and $ L/(T \oplus S) \simeq T^{*}/T$, $\frac{al + bt}{2n}$ is mapped to $\frac{al}{2n}$ and $\frac{bt}{2n}$. Therefore, $(a, 2n) = (b, 2n) = 1$. So, we may write from the first that $L/(T \oplus S) = \langle \frac{l + bt}{2n} \rangle$, 
where $(b, 2n) = 1$. Note that there is a natural one to one correspondence between over lattices of $S \oplus T$ and the subgroups of $(S^{*} \oplus T^{*})/(S \oplus T) = S^{*}/S \oplus T^{*}/T$. Therefore 
$b\, \text{mod}\, 2n$ is uniquely determined by $L$.   
\proclaim{Lemma (4.3)} Let $Y \in \Cal S_{T}(X)$.
\roster 
\item 
Assume that $n \not= 1$. Then $(f_{Y} \oplus g_{Y})(H^{2}(Y, \bold Z)) \not= (f_{Y} \oplus (-g_{Y}))(H^{2}(Y, \bold Z))$ and there exists unique $i = i(Y) \in \{1, \cdots , M \}$ such that $\tau_{X_{i}}(H^{2}(X_{i}, \bold Z)) = (f_{Y} \oplus (-g_{Y}))(H^{2}(Y, \bold Z))$. 
\item 
Assume that $n = 1$. Then $(f_{Y} \oplus g_{Y})(H^{2}(Y, \bold Z)) = (f_{Y} \oplus (-g_{Y}))(H^{2}(Y, \bold Z))$. 
\endroster 
\endproclaim 
\demo{Proof} We can write $H^{2}(Y, \bold Z)/(\text{NS}(Y) \oplus T_{Y}) = \langle \frac{l_{Y} + bt_{Y}}{2n} \rangle$, where 
$l_{Y}$ is the same as before 
and $t_{Y}$ is an element of $T_{Y}$ such that $g_{X}(t_{Y}) = t$ and $b$ is an integer such that $(b, 2n) = 1$. Then, $(f_{Y} \oplus g_{Y})(H^{2}(Y, \bold Z))/(S \oplus T) = \langle \frac{l + bt}{2n} \rangle$ and 
$(f_{Y} \oplus -g_{Y})(H^{2}(Y, \bold Z))/(S \oplus T) = \langle  \frac{l - bt}{2n} \rangle$. Then, $(f_{Y} \oplus g_{Y})(H^{2}(Y, \bold Z)) = (f_{Y} \oplus -g_{Y})(H^{2}(Y, \bold Z))$ in $S^{*} \oplus T^{*}$ if and only if $2b \equiv 0\, \text{mod}\, 2n$. Since $(b, 2n) = 1$, this implies $2 \equiv 0 \text{mod}\, 2n$. However, this is possible only when $n = 1$. Thus, we get the first part of (1). The last part of (1) is nothing but the definition of the set $\{X_{j}\}_{j=1}^{M}$. If $n=1$, one has $b =  - b\, \text{mod}\, 2$ and the two lattices are the same in $S^{*} \oplus T^{*}$. \qed \enddemo
\proclaim{Lemma (4.4)} 
\roster 
\item 
For each $j \in \{1, \cdots , M\}$, there is exactly one $i \in \{1, \cdots , M \}$ such that $i \not= j$ and such that $X_{i} \simeq X_{j}$ when $n \not= 1$. When $n =1$, $X_{i} \simeq X_{j}$ if and only if $i = j$.  
\item
In particular, $\vert \Cal S_{T}(X) \vert = M/2$ if $n \not= 1$ and $\vert \Cal S_{T}(X) \vert = M$ if $n = 1$. 
\endroster
\endproclaim 
\demo{Proof} (2) follows from (1) and (4.2). Let us show (1). If $n \not= 1$, then setting $\rho_{j} := (f_{X_{j}} \oplus -g_{X_{j}}) \vert H^{2}(X_{j}, \bold Z)$, one has 
$\rho_{j}(H^{2}(X_{j}, \bold Z)) = L_{i(j)}$ for some $i(j) \not= j$ by (4.3). For this $i = i(j)$, the map $\tau_{i}^{-1} \circ \rho_{j} : H^{2}(X_{j}, \bold Z) \rightarrow H^{2}(X_{i}, \bold Z)$ is a Hodge isometry. Therefore $X_{i} \simeq X_{j}$ by the global Torelli Theorem. 
\par 
\vskip 4pt 
Assume that $X_{k} \simeq X_{j}$. Let $\varphi : X_{k} \rightarrow X_{j}$ be 
an isomorphism. Note that $l_{X_{k}}$ is the unique ample generator of $\text{NS}(X_{k})$. 
Then one has $\varphi^{*}(l_{X_{j}}) = l_{X_{k}}$, $\varphi^{*}(T_{X_{j}}) = T_{X_{k}}$ and $\varphi^{*}(H^{2}(X_{j}, \bold Z)) = H^{2}(X_{k}, \bold Z)$.  
Set $\rho := \tau_{X_{k}} \circ \varphi^{*}$ This $\rho$ is a Hodge isometry and satisfies
$$\rho(l_{X_{j}}) = \tau_{X_{k}}(l_{k}) = l\,\, \text{and}\,\, \rho(H^{2}(X_{j}, \bold Z)) = \tau_{X_{k}}(H^{2}(X_{k}, \bold Z)) = L_{k}\, .$$ 
On the other hand, by (4.1), we have exactly two Hodge isometries from $H^{2}(X_{j}, \bold Z)$ into $S^{*} \oplus T^{*}$ which mapps $l_{X_{j}}$ to $l$. They are $(f_{X_{j}}, g_{X_{j}})$ and $(f_{X_{j}}, -g_{X_{j}})$. Thus, 
$\rho = (f_{X_{j}}, g_{X_{j}})$ or $\rho = (f_{X_{j}}, -g_{X_{j}})$. In the first case, we have $L_{k} = L_{j}$ and in the second case we have $L_{k} = L_{i(j)}$. By the definition of $L_{*}$, we have then $k = j$ and $k = i(j)$ 
respectively. Then, by the definition of the set $\{X_{i}\}_{i=1}^{M}$, 
we have $X_{k} = X_{j}$ and $X_{k} = X_{i(j)}$ respectively. This completes 
the proof of (1) in the case $n \not= 1$. Proof for the case $n=1$ is exactly the same. \qed \enddemo
In order to complete the proof of (1.10), it now suffices to prove the following Lemma:
\proclaim{Lemma (4.5)} $M = (\vert ((\bold Z/4n)^{\times})_{2} \vert)/2$. \endproclaim
\demo{Proof} Recall that $M$ is the cardinality of the set consisting of the even unimodular 
over lattices $L$ of $S \oplus T$ such that $S$ and $T$ are both primitive in $L$. We write this set by $\Cal M$.
\par 
\vskip 4pt
Let $L \in \Cal M$. Then, as remarked before (4.3), taking an element $t$ such that $T^{*}/T = \langle \frac{t}{2n} \rangle$, we can write $L/(T \oplus S) = \langle \frac{l + bt}{2n} \rangle$. Then $\overline{b} := b\, \text{mod}\, 2n$ is an element of $(\bold Z/2n)^{\times}$ and is uniquely determined 
by $L$. 
\par
\vskip 4pt
Conversely, if an over lattice $L$ of $S \oplus T$ satisfies $L/(S \oplus T) = \langle \frac{l + bt}{2n} \rangle$ for some natural number $b$ such that $(b, 2n) = 1$, i.e. $\overline{b} \in (\bold Z/2n)^{\times}$, then $\vert \text{det}\, L \vert = 1$ and both $S$ and $T$ are primitive in $L$. 
Therefore, there is one-to-one correspondence between $\Cal M$ and the set of element $\overline{b}$ of $(\bold Z/2n)^{\times}$ such that the overlattice $S + T + \langle \frac{l + bt}{2n} \rangle$ is even and integral. We write the latter set by $\Cal D$. Then $\vert \Cal M \vert = \vert \Cal D \vert$ and we may calculate $\vert \Cal D \vert$. Since $(\frac{t}{2n}, u) \in \bold Z$ and $(\frac{l}{2n}, u) \in \bold Z$ if $u \in S \oplus T$ (by the definition of the dual) and since $S$, $T$ are both even, we see that $L = S + T + \langle \frac{l + bt}{2n} \rangle$ is even and integral if and only if $(\frac{l + bt}{2n})^{2}$ is an even number. Since $(\frac{t}{2n}, t) \in \bold Z$, we have $(t^{2}) = 2nc$ for some integer $c$. Then $b$ satisfies that 
$$(\frac{l + bt}{2n})^{2} = \frac{2n + b^{2}(t^{2})}{4n^{2}} 
= \frac{1 + b^{2}c}{2n}\, .$$
Thus, the condition $b\, \text{mod}\, 2n \in \Cal D$ is 
equivalent to
$$1 + b^{2}c = 0\, \text{mod}\, 4n\,\, ,\text{i.e.}\,\,  b^{2}c = -1\, \text{in}\, \bold Z/4n.$$ 
There is at least one $b\, \text{mod}\, 2n \in \Cal D$ corresponding to the original $X$. We write this $b$ by $b_{0}$. Then $b_{0}^{2}c = -1\,  \text{in}\, \bold Z/4n$. In particular, both $b_{0}$ and $c$ are unit elements in $\bold Z/4n$. Therefore the condition $b\, \text{mod}\, 2n \in \Cal D$ is equivalent to the condition $(\frac{b}{b_{0}})^{2} = 1$ in $\bold Z/4n$, that is, $b$ is written as $b_{0}u$ for an integer $u$ which gives an element of 
$(\bold Z/4n)^{\times}$ of order at most $2$. Since $b_{0}$ is invertible in $\bold Z/4n$, we see that $b_{0}$ is invertible in $\bold Z/2n$ as well. Then $b_{0}u = b_{0}v$ in $\bold Z/2n$ if and only if $u = v$ in $\bold Z/2n$, i.e. $v = u$ or $v = u + 2n$ in $\bold Z/4n$. Note that $u + 2n \not= u$ in $\bold Z/4n$. Note also that if 
$u$ is a unit element in $\bold Z/4n$, then $u + 2n$ is also a unit element in $\bold Z/4n$. Indeed, if $uc = 1$ in $\bold Z/4n$, then $u = c = 1$ in $\bold Z/2$. Therefore 
$$(u + 2n)(c + 2n) = uc + 2n(u+c) + 4n^{2} = 1\, \text{in}\, \bold Z/4n\, .$$ 
Hence $M = \vert \Cal D \vert = \vert \{b_{0}u\, \text{mod}\, 2n \}\vert 
= \vert ((\bold Z/4n)^{\times})_{2} \vert/2$. \qed \enddemo 
Now combining (4.4) and (4.3), we obtain (1.10). The explicit formula 
follows from the Chinese remainder Theorem and the fact $(\bold Z/2^{2})^{\times} \simeq \bold Z/2$, $(\bold Z/2^{a})^{\times} \simeq \bold Z/2 \oplus \bold Z/2^{a-2}$ if $a \geq 3$ and $(\bold Z/p^{e})^{\times} \simeq \bold Z/p^{e-1}(p-1)$ if $p \geq 3$ is a prime. \qed 
\head
{Appendix} 
\endhead 
In this appendix, we shall remark a generalization of (1.7) for hyperk\"ahler manifold of higher dimension.
\par 
\vskip 4pt 
By a hyperk\"ahler manifold, we mean a simply-connected compact K\"ahler 
manifold $X$ which admits a everywhere non-degenerate holomorphic 
$2$-form $\omega_{X}$ and satisfies $H^{2, 0}(X) = \bold C \omega_{X}$. Let 
$X$ be a hyperk\"ahler manifold. $X$ is of even dimensional. By Beauville [Be], $X$ admits a non-degenerate primitive integral symmetric bilinear form $(*, *) : H^{2}(X, \bold Z) \times H^{2}(X, \bold Z) \rightarrow \bold Z$ of signature $(3, B_{2}(X) - 3)$. By his construction, this bilinear form together with the Hodge decomposition $H^{2}(X, \bold C) = H^{1, 1}(X) \oplus \bold C \omega_{X} \oplus \bold C \overline{\omega_{X}}$ forms a weight two Hodge structure on $H^{2}(X, \bold Z)$. We define the transcendental lattice $T_{X}$ of $X$ to be the orthogonal lattice of $\text{NS}(X) \simeq \text{Pic}\, X$ in $H^{2}(X, \bold Z)$ with respect to this bilinear form. $T_{X}$ is the minimal primitive sub-lattice of $H^{2}(X, \bold Z)$ whose $\bold C$-linear extension contains $\omega_{X}$. As is remarked in [Hu1] (all of which results are now avaiable by [Hu2]), $H^{2}(X, \bold Z)$ and $H^{2}(X', \bold Z)$ are Hodge isometric if $X$ and $X'$ are birational. 
\par 
\vskip 4pt 
If $X$ is the Hilbert scheme $\text{Hilb}^{n}(S)$ of $0$-dimensional subschemes of length $n (\geq 2)$ on a K3 surface $S$ or its (K\"ahler) deformation, 
then $X$ is a hyperk\"ahler manifold of dimension $2n$ and satisfies 
$(H^{2}(X, \bold Z), (*, *)) \simeq \Lambda \oplus \langle -2(n -1) \rangle$. Here $\Lambda$ is a K3 lattice. 
\par
\vskip 4pt
\proclaim{Proposition A} Let $n$ be any natural number. For any given natural number $N$, there are $N$ 
mutually non-birational projective hyperk\"ahler manifolds 
$\{X_{j}\}_{j=1}^{N}$ of dimension $2n$ such that $(T_{X_{j}}, \bold C \omega_{X_{j}})$ are mutually Hodge isometric but $\text{NS}(X_{j})$ are not isomorphic to one another. 
\endproclaim 
\demo{Proof} Let $S_{i}$ and $T_{i}$ ($i = 1, \cdots , N$) be the lattices found in (2.9) and (2.12). 
Consider the sublattices $\tilde{T}_{i} := T_{i} \oplus \langle -2(n -1) \rangle \subset \Lambda \oplus \langle -2(n -1) \rangle$. Then, $\tilde{T}_{i}$ are primitive, of signature $(2, 19)$ and are isomorphic to one another by the construction and (2.11). (We remark that $l(A_{\tilde{T}_{i}}) \leq l(A_{T_{i}}) + 1 \leq 3$ so that we can apply (2.11) to see these $\tilde{T}_{i}$ are isomorphic.) Fix a lattice $\tilde{T}$ isomorphic to all $\tilde{T}_{i}$ and fix primitive embeddings $\tilde{\varphi}_{i} : \tilde{T} \simeq \tilde{T}_{i} \subset \Lambda \oplus \langle -2(n -1) \rangle$. As in (2. 13), we also take a positive definite $2$ dimensional subspace $\tilde{P}$ of $\tilde{T} \otimes_{\bold Z} \bold R$ 
such that $\tilde{P}^{\perp} \cap \tilde{T} = \{0\}$. One can then obtain a weight $2$ Hodge structure on $\tilde{T}$, say, $(\tilde{T}, \bold C \tilde{\omega})$ as in (2.14). Then, by applying the surjectivity of the period mapping by Huybrechts [Hu1] for the weight two Hodge structures $(\Lambda \oplus \langle -2(n -1) \rangle, \bold C \tilde{\varphi}_{i}(\tilde{\omega}))$ and repeating exactly the same argument as in (2.14), one can find for each $i$ a marked hyperk\"hler manifold $(X_{i}, \tau_{i})$ (equivalent to $\text{Hilb}^{n}(S)$ under deformation) such that $(T_{X_{i}}, \bold C \omega_{X_{i}}) \simeq (\tilde{T}_{i}, \bold C \tilde{\varphi}_{i}(\tilde{\omega}))$ 
and $\text{NS}(X_{i}) \simeq S_{i}$ via $\tau_{i}$. Since $S_{i}$ 
is of signature $(1,1)$, $X_{i}$ is projective by the projectivity criterion 
due to Huybrechts [Hu1] (See [Hu2] for a correction of proof). Therefore, 
these $X_{i}$ give the result. \qed \enddemo 
\proclaim{Remark} 
\roster 
\item In [Yo], Yoshioka finds K3 surfaces $X \not\simeq Y$ such that 
$\text{Hilb}^{2}(X) \simeq \text{Hilb}^{2}(Y)$. In addition, there are several lattices $S$, $S'$ and $M$ such that $S \not\simeq S'$ but 
$S \oplus M \simeq S' \oplus M$. So, Hilbert schemes of K3 surfaces found in (1.7) may not satisfy the condition of Propsition A. 
\item Namikawa [Nm] finds a counter example of the birational 
injectivity of the period mapping for hyperk\"ahler manifolds of dimension 4. Therefore, the argument for the finiteness (1.9) in Section 3 cannot be exploited in higher dimensional hyperk\"ahler manifolds. \qed
\endroster 
\endproclaim


\Refs
\widestnumber\key{BPV}
\ref 
\key Be
\by A. Beauville
\paper Vari\'et\'es K\"ahleriannes dont la premi\'ere classe de Chern est 
nulle 
\jour J. Diff. Geom. 
\vol 18
\yr 1983
\pages 755--782
\endref

\ref
\key BO 
\by  A. Bondal, D. Orlov
\paper Reconstruction of a variety from the derived category and groups of autoequivalences, math.AG/9712029
\endref

\ref 
\key BM1
\by T. Bridgeland, A. Maciocia
\paper Complex surfaces with equivalent derived categories
\jour Math. Z. 
\vol 236
\yr 2001
\pages 677--697
\endref

\ref
\key BM2 
\by  T. Bridgeland, A. Maciocia
\paper Fourier-Mukai transforms for K3 and elliptic fibrations, 
math.AG/9908022
\endref

\ref 
\key BPV
\by W. Barth, C. Peters, A. Van de Ven
\paper Compact complex surfaces 
\jour  Springer-Verlag 
\yr 1984
\endref

\ref 
\key Cl
\by A. Caldararu
\paper Non-fine moduli spaces of sheaves on K3 surfaces, math.AG/0108180
\endref

\ref 
\key Cs
\by J. W. S. Cassels
\paper Rational quadratic forms
\jour Acdemic Press
\yr 1978
\endref

\ref
\key GM
\by S.I. Gelfand, Y.I. Manin
\paper Methods of homological algebra
\jour Springer-Verlag
\yr 1991
\endref

\ref 
\key Hu1
\by D. Huybrechts
\paper Compact Hyperk\"ahler manifolds: Basic results 
\jour Invent. math. 
\vol 135
\yr 1999
\pages 63--113
\endref

\ref 
\key Hu2
\by D. Huybrechts
\paper Erratum to the paper:Compact Hyperk\"ahler manifolds: Basic results, 
math.AG/0106014 
\endref

\ref
\key Iw
\by H. Iwaniec
\paper Almost-primes represented by quadratic polynomials
\jour Invent. math.
\vol 47
\yr 1978
\pages 171--188
\endref

\ref
\key Mu
\by S. Mukai
\paper On the moduli space of bundles on K3 surfaces I, in: Vector bundles on algebraic varieties
\jour Oxford Univ. Press
\yr 1987
\pages 341--413
\endref

\ref
\key Nm
\by Yo. Namikawa
\paper Counter-example to global Torelli problem for irreducible symplectic manifolds, math.AG/0110114
\endref

\ref
\key Nr
\by W. Narkeiwicz
\paper Elementary and anlytic theory of algebraic numbers, 2nd edition
\jour Springer-Verlag
\yr 1990
\endref

\ref 
\key Ni
\by V. V. Nikulin
\paper Integral symmetric bilinear forms and some of their geometric 
applications
\jour Math. USSR Izv. 
\vol 14
\yr 1980
\pages 103--167
\endref

\ref 
\key Or
\by D. Orlov
\paper Equivalences of derived categories and K3 surfaces, math.AG/9606006
\endref

\ref 
\key PS
\by I. Piatetski-Shapiro and I. R. Shafarevich 
\paper A Torelli Theorem for algebraic surfaces of type K3
\jour Math. USSR Izv.
\vol 5
\yr 1971
\pages 547--587
\endref

\ref 
\key Th
\by R. Thomas 
\paper A holomorphic Casson invriant for Calabi-Yau 3-folds and bundles on K3 fibrations, math.AG/9806111
\endref  

\ref 
\key Yo
\by K. Yoshioka 
\paper Moduli spaces of stable sheaves on abelian surfaces, math.AG/0009001
\endref
\endRefs
\enddocument